\documentclass{amsart}
\usepackage{amsmath,amssymb,amsfonts,latexsym}
\let\nd=\noindent
\let\la=\langle
\let\ra=\rangle
\let\wt=\widetilde
\def\hol#1{\ensuremath{\rm Hol \, } (#1)}
\def\inv{\ensuremath{\rm Inv}(\Gamma ,G,m)}
\def\epi{\ensuremath{\rm Epi}(\Gamma ,G)}
\def\rep{\ensuremath{\rm Rep}(G ,s)}

\def\affeq{\ensuremath{{E}_{\rm Aff}^{\rm cf}}}
\def\diffeq{\ensuremath{{E}_{\rm Diff}^{\rm cf}}}

\def\w{\ensuremath{\rm w}}

\def\card{\ensuremath{\rm card \,} }
\def\rr#1{{\mathbb R  }^ {\hskip0.03cm #1}}
\def\rrr{{\mathbb R }}
\def\zz#1{{\mathbb Z }_{#1}}

\def\zzz{{\mathbb Z}}
\def\ccc{{\mathbb C}}
\def\qqq{{\mathbb Q}}

\def\ff{{\mathbb F}}
\def\tors{\ensuremath{\rm Tors \,  }}

\def\diff{\ensuremath{\rm Diff}\, }
\def\aff#1{\ensuremath{\rm Aff}\,(#1)}
\def\sl{\ensuremath{\rm SL}}
\def\gl{\ensuremath{\rm GL}}
\def\flbx{FLB_{\ccc }(X)}
\newtheorem{theorem}{Theorem}[section]
\newtheorem{proposition}{Proposition}[section]
\newtheorem{lemma}{Lemma}[section]
\newtheorem{E}{Example}[section]

\newtheorem{corollary}{Corollary}[section]
\newtheorem{D}{Definition}[section]
\newenvironment{defn}{\begin{D} \rm }{\end{D}}
\newtheorem{R}{Remark}[section]

\begin{document}

\bibliographystyle{plain}

%%%%%%%%%%%%%%%%%%%%%%%%%%%%%%%%%%%%%%%%%%%%%%%%%%
\title{Topological and affine structure of complete flat manifolds}
\author{Micha\l \ Sadowski}
\email{msa@delta.math.univ.gda.pl}
\address{   
 Department of Mathematics 
Gda\'nsk University 
80-952 Gda\'nsk, Wita Stwosza 57
Poland   } 

\date{\today}

\begin{abstract} 
The results of the paper concern the topological structure of 
complete riemannian manifolds with cyclic holonomy groups 
and low-dimensional orientable complete flat manifolds. 
We also discuss related results such as 
the affine classification of orientable complete flat $4$-manifolds, 
an algebraic criterion of an affine equivalence and 
the relationship between holonomy homomorphisms and 
some algebraic and geometric invariants. 
\end{abstract} 

\maketitle

\section{Introduction} 
The aim of this paper is to collect some results concering 
topological and affine structure of 
complete flat manifolds({\it cf-manifolds}). 
We recall particularly important known results in Section  
\ref{sec-holonomy-repr}. 
The others seem to be new. 
Complete flat manifolds play a particular role in geometry. 
On one hand they are natural generalizations of euclidean spaces, 
having the same local properties. 
On the other hand the study of some complicated questions, 
arising in differential geometry and related fields, 
often starts with the examination of them in the case 
of manifolds of constant curvature. 
%%%%%%%%% 
 A cf-manifold  is  the orbit space 
$M=\wt M/\Gamma $ 
of a properly discontinuous, isometric, and free  action 
of a discrete group $\Gamma $ on an euclidean space 
$\wt M. $ 
The  {\it holonomy homomorphism} of $M$ is the
map $\Phi $ carrying $\gamma \in
\Gamma $ onto linear part of $\gamma $ and 
$\Phi (\Gamma )$ is the holonomy group of $M.$ 
The linear isometry $\Phi (\gamma )$ can be written as 
$ \Phi _X (\gamma ) \times \Phi _U (\gamma ),$ 
where $\Phi _X(\gamma )$ acts on the universal covering space 
$\wt X$ of the totally geodesic submanifold $X$ of $M$, 
homotopy equivalent to $M,$ 
and $\Phi _U(\gamma )$ acts on the orthogonal complement 
of $\wt X.$ 
If $M$ is compact, then $M=X$ and $\Phi _X=\Phi .$  
In this case we say that $M$ is a 
{\it Bieberbach manifold}. 
If $M$ is noncompact, then $M$  
is determined by the Bieberbach group 
$\Gamma $ and the vertical holonomy homomorphism 
$\Phi _U$ so that the theory of cf-manifolds can be treated as 
the theory of orthogonal representations of Bieberbach groups. 
The manifold $M$ is the total space of a 
flat riemannian vector bundle $\xi :M\to X$ 
whose structure group can be reduced to $\Phi _U(\Gamma ) . $
The main
difference between the noncompact case and the compact one is that
$\Phi (\Gamma )$ is not always finite and $\Phi (\Gamma )$ is not
a topological invariant of $M$ (cf. Section \ref{sec:affeq} below).

\vskip0.2cm 
The main difficulty in the classification of Bieberbach  manifolds 
is that it is based on the classification of conjugacy classes 
of finite subgroups of $\gl (n,\zzz ).$ 
This is a hard problem in integral representaion group theory, 
solved for cyclic groups of prime order $p$, 
for cyclic groups of order $p^2$, and for small values of $n$ 
only. 
Bieberbach manifolds are described when the holonomy groups of them are cyclic
groups of prime order and when the dimensions of them are smaller than 
$7.$ 
The complexity of the problem is also connected with the fact that 
the number 
$\nu _f (n)$ of $n$-dimensional closed flat  manifolds
increases rapidly with $n.$ It is known that $\nu _f (2)=2,$ $\nu
_f(3) =10$ (\cite[Section 3.5]{wolf}), $\nu _f(4) =74,$ $\nu _f(5)
=1060$ (\cite{cidschulz}), and $\nu _f(6) =38746$
(\cite{cidschulz}).

\vskip0.2cm 
In this paper we study cf-manifolds in two particular important cases: 
when they have cyclic holonomy groups and when their dimensions are 
smaller than $5.$ 
These particular cases are starting points in the investigation of more 
general ones. 
In the affine classification it is convenient to use a 
criterion of affine equivalence formulated in Section \ref{sec:affeq}. 
%% parallel spinors on cf-manifolds and we describe how to calculate some   
%% 
We also answer some natural questions concerning 
algebraic and geometric invariants that are used in the paper. 
Flat manifolds were investigated in many books and papers. 
Only few of them deal with the noncompact case 
(see e.g. \cite{gw}, \cite{wolf}, \cite{ktflat}, \cite{ow}, 
\cite{tapp-thesis}, \cite{tapp}, \cite{wilking}). 
Related results can be found in papers dealing with flat 
vector bundles, for instance in 
\cite{cg}, \cite{ds}, \cite{gh}, \cite{gw}, \cite{ht}, and \cite{kt}. 

\vskip0.3cm
Throughout this paper the following notation will be used.
The universal covering space of a topological space $Y$
will be denoted by $\widetilde Y.$
If $H$ is a group, $h_1,...,h_k\in H$, then $\langle h_1,...,h_k\rangle $
will denote the subgroup of $H$ generated by $h_1,...,h_k$. 
$\diff (M)$ is the group of diffeomorphisms
of a smooth manifold $M$
and  $\aff V$ is the group of affine diffeomorphisms of $V.$
By $\hol  X$ we denote the holonomy group of
a riemannian manifold $X$.
The symbol $ \Pi _{ab} $ stands for projection of a Bieberbach group 
$\Gamma $ onto its abelianization $H_1(\Gamma )$  
and $E\xi $ for the total space of a vector bundle $\xi .$ 
Real and complex $d$-dimensional trivial vector bundles will be denoted
by $\Theta ^d _\rrr $ and
$\Theta ^d _\ccc  $ respectively.

%%%%%%%%%%%%%
\section[Holonomy representations and vector bundles]{Holonomy representations
and vector bundles}
\label{sec-holonomy-repr}

%%% \vskip0.7cm
The aim of the first part of this section is to 
reformulate some known results in a more
convenient for our purposes form.
The results will be expressed in terms of vector bundles associated
with coverings.
To describe these bundles consider a covering map
$\Pi : \widehat X \to X $ and its
covering transformation group $H$.
Let  $\ff $ be the field of reals or the field of complex numbers,
and  let $\rho :H \to \gl (s,\ff ) $ be an $s$-dimensional
representation of $H.$
Take the diagonal action $h(x,u)=(hx,\rho (h)u)$ on
$\widehat X \times \ff ^s$ and the orbit space $\Pi [\rho ]=
(\widehat X\times \ff ^s)/H.$
Let $p: \Pi [\rho ] \to X $ be the map
determined by the projection $p_0 :
\widehat X \times \ff ^s \ni (x,u) \to x \in \widehat X $.
Then the triple $(\Pi [\rho
], X,p) $ is a vector bundle associated to the principal bundle
$\Pi $ with typical fiber $\ff ^s.$ 
In the sequel we identify
this bundle with $\Pi [\rho ]. $
We  often use the following.

\begin{proposition}\label{propclosedflatsubmanifold}
Let $M$ be a connected $m$-dimensional cf-manifold
and
let $X$ be a  closed totally geodesic submanifold of $M$,
homotopy equivalent to $M.$
Denote  $\Phi (\Gamma )$ by $H$, $\dim X $ by $n$,
and $m-n$ by $s.$
Then there are
a riemannian covering $\Pi : \widehat X \to X$
and
an orthogonal representation $\rho :H\to {\rm O}(s)$
such that
$M=\Pi [\rho ].$
The projection
$p: \Pi [\rho ] \to  X$ is affine.
The action of $\gamma \in \Gamma $ on
$\widetilde M$ can be written as
$(t_{v(\gamma )} \circ \Phi _X (\gamma ) )\times \Phi _U (\gamma ),$
where $\Phi _X (\gamma ) :\widetilde X \to \widetilde X $
and
$\Phi _U (\gamma ) :\widetilde X ^\perp \to \widetilde X ^\perp .$
\end{proposition}

There are different proofs of Proposition \ref{propclosedflatsubmanifold}. 
For instance it follows from the arguments used in the proof of 
Theorem 3.3.3 in \cite[ch. 3]{wolf}. 
The homomorphism $\Phi _U$ is the {\it vertical holonomy homomorphism} 
od $M.$

\begin{corollary}\label{corsplitcommutativeflatbundles}
If the holonomy group $H$ of $M$ is abelian, 
then the bundle
$p: M \to X$ is isomorphic to
$$
\bigoplus _{j} \Pi [\rho _j]
\oplus \bigoplus _{k} \Pi [\lambda _k] ,$$
where $\rho _j$ are 1-dimensional complex representations
of $H$ and
$\lambda _k$
are 1-dimensional real representations of $H.$
\end{corollary}

%%%%%%%%%%%%%%%%%%%%%%%%%%%%%%
Proposition \ref{propclosedflatsubmanifold} and the following observation 
show that cf-manifolds correspond to orthogonal representations 
of Bieberbach groups. 

\begin{proposition}\label{prop-cf-man-and-repr} 
If $X=\widetilde X/\Gamma $ is a Bieberbach manifold, 
$\rho :\Gamma \to {\rm O}(s)$ is an orthogonal representation of $\Gamma ,$  
and  
the action of $\gamma \in \Gamma $ on 
$\wt \Gamma \times \rrr ^s$ is given by the formula 
$$ \gamma (x,u) =(\gamma x,\rho (\gamma )u ) $$ 
then the orbit space 
$ M(\rho ) = (\widetilde X \times \rrr ^s)/\Gamma $ 
is a complete flat manifold homotopy equivalent to $X.$ 
\end{proposition} 

If $\rho ^\prime :\Gamma \to {\rm O}(s)$ is 
another orthogonal representation of $\Gamma ,$
then it is natural to ask when 
$M(\rho )$ is diffeomorphic (affinely diffeomorphic) to 
$M(\rho ^\prime )$. We consider this question in the next sections. 
In particular, we give an algebraic criterion of an affine equivalence. 

%%%%%%%%%%%%%%%%%%%%%%%%%%%%
\vskip0.2cm 
The following  three simple lemmas are known.

\begin{lemma}\label{lem-complexif-and-restr-for-bundles}
{\rm (cf. \cite[ch. 16, \S \ 11]{husemoller}).}
If $Y$ is a topological space,
 $\xi :M\to Y $ is a real vector bundle and
$\zeta : V\to Y $ is a complex  vector bundle, then
$$  (r \circ c)([\xi ]_\rrr ) =2[\xi ]_\rrr \quad \hbox {and} \quad
  (c \circ r)([\zeta ])=[\zeta \oplus \zeta ^\ast ] .$$
\end{lemma}
The maps $r$ and $c$ 
carrying 
$ [\zeta ]_\ccc \in KU(X) $ onto $[\zeta ]_\rrr \in KO(X) $ 
and $ [\xi ]_\rrr \in KO(X)$ onto $[\xi \otimes \ccc ]_\ccc \in KU(X)$ 
are called the restriction and the complexification.

\begin{lemma}\label{lem-stability-for-bundles}
{\rm (cf. \cite[ch. 8, Theorem 2.6]{husemoller}).}
Let $Y$ be a topological space homotopy equivalent to a finite
$n$-dimensional CW-complex.

\vskip0.1cm
\nd
{\bf a)}
If  $\xi $  and $\xi ^\prime  $ are two real vector bundles
over $Y$,   $\dim _\rrr \xi =\dim _\rrr \xi ^\prime \geq n+1,$
and $[\xi ]_\rrr =[\xi ^\prime ]_\rrr $, then
$\xi \cong \xi ^\prime .$

\vskip0.1cm
\nd
{\bf b)}
If  $\zeta  $  and $\zeta ^\prime  $ are two complex vector bundles
over $Y$,   $\dim _\rrr \zeta =\dim _\rrr \zeta ^\prime \geq n,$
and $[\zeta  ]_\ccc  =[\zeta ^\prime ]_\ccc  $, then
$\zeta \cong \zeta ^\prime .$
\end{lemma}

%%%%%%%%%%%%%%%%%%%%%

\begin{lemma}\label{lem-flat-bun-are-torsion-elem-in-k-theory}
Let $X^n$ be a closed flat manifold.
Assume that $\xi :M\to X $ is an $s$-dimensional real flat vector bundle and
$\zeta : V\to X $ is a $q$-dimensional complex flat vector bundle.
Then

\vskip0.1cm
\nd
{\bf a)}  $M$ and $V$ have the structures of complete flat manifolds,

\vskip0.1cm
\nd
{\bf b)}  the image of the 
total Chern class ${\rm c}(\zeta )$ in $H^\ast (X, \qqq )$
is equal to $1,$

\vskip0.1cm
\nd
{\bf c)}  $ch(\zeta )=1,$

\vskip0.1cm
\nd
{\bf d)}  $[\zeta ] -[\Theta _\ccc ^q] \in \tors KU(X),$

\vskip0.1cm
\nd
{\bf e)}  $[\xi ] -[\Theta _\rrr ^s] \in \tors KO(X).$
\end{lemma}

%%%%%%%%%%%%%%%%%%%%%%%%%%%%%%%%%%%%%%%
The proof of a) can be found in \cite{milnor} (see also \cite{ow}).
Parts b) and c) are consequences of the flatness of $\zeta $,
which implies that Chern forms representing Chern classes of
$\zeta $ are equal to  $0$ (see e.g.
\cite[ch. 11]{kn} and
\cite[Appendix C, Corollary 2]{milnorstasheff}).
Parts d) and e)  follows from the fact that 
$ch : KU(X) \otimes \qqq \to H^{ev} (X,\qqq )$ is a
monomorphism (cf. \cite[ch. 5, Theorem 3.25]{karoubi},
\cite[ch. 1, Section B, \S \  4]{dyer}), from the equality 
 $ch([\zeta ] - [\Theta _\ccc ^q])=0$, and from Lemma 
\ref{lem-complexif-and-restr-for-bundles}.  
%%%%%%%%%%%%%%%%%%%%%%%%%%%%%%%%%%%%%%%%%%%%

\vskip0.2cm 
Now we formulate a particular case of a result of Wilking 
(\cite[Corollaries 6.4 and 6.5]{wilking}). 
It is a generalization of the second Bieberbach theorem.

\begin{theorem}\label{thm-finiteness-for-cf-man} 
The set  $\diffeq (m)$ of  diffeomorphism classes of 
\hyphenation{dimen-sional} 
$m$-dimensional 
cf-manifolds is finite and each element of 
$\diffeq (m)$ contains a cf-manifold with finite holonomy group. 
\end{theorem} 

An affine variant of Theorem \ref{thm-finiteness-for-cf-man} 
is false because the holonomy groups of affinely diffeomorphic flat 
manifolds are isomorphic. 

%%%%%%%%%%%%%%%%%%%%%%%%%%%%%%%%%
\section{Complete flat manifolds with cyclic holonomy groups}
\label{sec-cfch-manifolds}

The aim of this section is to give a topological
classification of complete flat manifolds whose holonomy groups
are cyclic and whose fundamental groups are isomorphic to a fixed
group $\Gamma .$ 

\begin{defn}
A flat riemannian manifold with cyclic holonomy group
will be called a {\it fch-manifold}.
A {\it nolcyc bundle} is a nonorientable line bundle
$L :EL \to X$ such that $EL$ is a complete fch-manifold
and $X$ is a totally geodesic submanifold of $ EL $.
\end{defn}
The main here are the following.  

\begin{theorem}\label{thmflatvectorbundlescyclicholonomy}
Let $M$, $X$, $\Pi [\rho ],$ $\Gamma ,$ and $s$ be as in
Proposition \ref{propclosedflatsubmanifold} and let 
$\epsilon (\gamma ) =\det \Phi _U (\gamma ).$ 
Assume that $\Phi (\Gamma )$ is a cyclic group.

\vskip0.1cm
\nd
{\bf a)}
If $\, \Pi [\rho ]$ is an orientable bundle, then $\Pi [\rho ] \cong 
\Theta _\rrr ^s.$

\vskip0.1cm
\nd
{\bf b)}
If $\, \Pi [\rho ]$ is a nonorientable bundle, then
$\Pi [\rho ] \cong \Theta _\rrr ^{s-1} \oplus \Pi [\epsilon ].$
\end{theorem}

\begin{theorem}\label{thmtopstructurefchmanifolds}
Let $M$ be a  complete, 
connected,  $m$-dimensional riemannian\hfill  manifold\hfill  with
\hfill  cyclic\hfill  holonomy\hfill group\hfill  and\hfill  let
\hfill  $X$\hfill  be\hfill  a\hfill  closed\\
 $n$-dimensional flat manifold
homotopy equivalent to $M$. Suppose that $L$ is a nolcyc-bundle
over $X$ and $\dim M > \dim X$.
Then $M$ is diffeomorphic either to
$X\times \rr {m-n}$
or to the total space of $(X\times  \rr {m-n-1} ) \oplus L.$
\end{theorem}

Theorem \ref{thmtopstructurefchmanifolds} shows that there are 
exactly two diffeomorphism classes of complete, noncompact fch-manifolds 
having the same fundamental group. 
It reduces the classification of complete fch-manifolds to 
the classification of compact ones. 
Recall that the classification of Bieberbach manifolds 
with a fixed cyclic holonomy group $C$ is a difficult problem 
solved only when $C$ has prime order (\cite{charlap} \cite{charlapb}). 

%%%%%%%%%%%
%%%%%%%%%%%%%%%%%%%%%%%%%%%%%%%%%%%%%%%%%%%%%%%%%%%%%%%%%%%%%%%%%%%%%%
\vskip0.2cm 
Theorem \ref{thmflatvectorbundlescyclicholonomy} 
is a consequence of the following algebraic property of 
the deck groups of fch-manifolds. 

\begin{lemma}\label{lemholonomyiszeroontorsion}
If $\Phi (\Gamma ) $ is a finite cyclic group, then
there are $a\in H_1(\Gamma )$ and a subgroup $B$ of
$ H_1(\Gamma )$
such that $H_1(\Gamma  )=\la a \ra \oplus B,$
$\la a \ra \cong \zzz ,$
$\Psi (\la a \ra ) = \Phi (\Gamma ) ,$
and ${\rm Tors}\, H_1(\Gamma ) \subset B\subset {\rm ker \, } \Psi . $
\end{lemma}

Theorem \ref{thmflatvectorbundlescyclicholonomy} 
implies that if the bundle $M\to X$ is orientable, then 
$M$ is diffeomorphic to $X\times \rr {m-n}$. 
The case when the bundle $M\to X$ is nonorientable is more difficult.  
It follows from the fact that any two nolcyc bundles 
$\lambda  $ and  $\lambda ^\prime $, with the same base space $X$,  
belong to the same orbit of the action of ${\rm Diff}(X)$ 
on $X.$ 

\vskip0.2cm 
Theorem \ref{thmflatvectorbundlescyclicholonomy} 
is a generalization of a result of Thorpe 
stating that  $X$ is parallelizable
or $TX$ is isomorphic to the direct sum 
of a trivial bundle and a line bundle (cf. \cite{th}). 
Using Theorem \ref{thmflatvectorbundlescyclicholonomy} 
it is easy to verify a more general version of the last statement. 

\begin{theorem}\label{thm-variant-of-thorpe-for-complete-flat}
Let $M$ be a complete fch-manifold.

\vskip0.1cm
\nd {\bf a)} If $M$ is orientable, then $M$ is parallelizable.

\vskip0.1cm
\nd {\bf b)} If $M$ is nonorientable, then
$TM \cong \Theta _\rrr ^{m-1} \oplus \lambda $ for some 
nolcyc-bundle  $\lambda $ over $M$.
\end{theorem}

\begin{corollary}
If $M$ is a complete flat manifold and $\hol {M}$ is a cyclic
group of odd order, then $M$ is parallelizable.
\end{corollary}

%%%%%%%%%%%%%%%%%%%%%%%
\section[Affinely equivalent complete flat manifolds]
{Affinely equivalent complete flat manifolds}
\label{sec:affeq} 

The aim of this section is to describe 
algebraic invariants corresponding to affine equivalence classes 
of noncompact cf-manifolds. 
For details and for the proofs we refer to  
\cite{msb}.

\begin{theorem}\label{thm-criterion-aff-eq}
Let $M$ and $M^\prime $ be two $m$-dimensional cf-manifolds
with isomorphic fundamental groups.
Let $X\subseteq M$ and $X^\prime \subseteq M^\prime $ be totally geodesic
submanifolds of $M$ and $M^\prime $ 
homotopy equivalent to $M$.
Assume that $n=\dim X,$ 
$\widetilde M=\widetilde X \times U  $ and 
$\widetilde M^\prime =\widetilde X^\prime  \times U^\prime . $
Let $\Phi _U : \pi _1(X) \to {\rm O}(m-n),$
$\Phi _U^\prime  :\pi _1(X^\prime )  \to {\rm O}(m-n)$
be the  vertical holonomy homomorphisms 
of $M$ and $M^\prime .$ 
Then the following conditions are equivalent.

\vspace{0.1cm}
\nd
{\bf (a)}
$M$ is affinely diffeomorphic to $M^\prime ,$

\vspace{0.1cm}
\nd
{\bf (b)}
there is an isomorphism $f: \pi _1(X) \to \pi _1(X^\prime )$
and a linear isomorphism $L:U\to U^\prime $ such that
$$ \Phi _U ^\prime (f(\gamma )) =
L \circ \Phi _U (\gamma ) \circ L^{-1} $$
for  $\gamma \in \pi _1(X)$. 
\end{theorem}

Let $m,$ $\Gamma ,$ and $n$ be as in the formulation of 
Theorem \ref{thm-criterion-aff-eq}. 
For a fixed discrete group $G$ consider the set 
${\mathcal I}(\Gamma ,G,m)$ of all pairs $(\epsilon ,\rho ),$ 
where $\epsilon :\Gamma \to G$ is an epimorphism and
$\rho :G \to {\rm O}(s)$ is a representation.

\begin{defn} 
Two elements $(\epsilon ,\rho )$ and $(\epsilon ^\prime ,\rho ^\prime )$
of ${\mathcal I}(\Gamma ,G,m)$  
are {\it equivalent} if there are $f\in {\rm Aut}(\Gamma ) $
and a linear isomorphism $L: \rr s \to \rr s $ such that
$$ L(\rho )\circ \epsilon =
\rho ^\prime  \circ \epsilon ^\prime \circ f ,$$ 
where $ L (\rho )(g) = L\circ \rho (g) \circ L^{-1}  .$ 
\end{defn} 

Let $\inv $ be the set of equivalence classes 
of the elements of ${\mathcal I}(\Gamma ,G,m) ,$ 
let $\epi $ be the set of
epimorphisms from $\Gamma $ to $G$, and let $\rep $ be the set of
conjugacy classes of representations of $G$ in $\rr s.$
Applying Theorem \ref{thm-criterion-aff-eq} we have. 

\begin{theorem}\label{thm-alg-inv}
If $\Gamma $ is  a Bieberbach group, then
there is a bijection $\nu : \inv \to \affeq (\Gamma ,G ,m).$
\end{theorem}

\begin{corollary}\label{cor-estimate-aff-eq-classes}
If $\vert G\vert <\infty ,$ then
$$\card  \affeq (\Gamma ,G ,m) \leq
\card \epi \, \card \rep
 .$$
\end{corollary}

Let $\affeq (\Gamma ,m)$ be the set of
affine diffeomorphism classes of
 $m$-dimensional complete flat manifolds
with the same fundamental group $\Gamma $
and let $n$ be as above.

\begin{theorem}\label{thm-uncountable-many-aff-eq}
If $m\geq n+2$ and $H_1(\Gamma ,\zzz )$ is infinite, then
$\affeq (\Gamma ,m)$ is uncountable.
\end{theorem}

%%%%%%%%%%%% 
\begin{proposition}\label{prop-infin-many-aff-classes} 
Let $\Gamma $ be a Bieberbach group. 
Then there are infinitely many affine equivalence classes of 
cf-manifolds whose fundamental groups are isomorphic to 
$\Gamma $ and whose holonomy groups are finite. 
\end{proposition}

%%%%%%%%%%%%%%%%%%%%
\section[Topological and affine classification]
{Topological and affine classification of low-dimensional\\
orientable cf-manifolds}
\label{sec:classif-low-dim}

The aim of this section is to describe topological  
and affine equivalence classes of 
cf-$4$-manifolds. For simplicity we deal with 
the orientable case. A nonorientable case is somewhat more 
complicated and will be considered elsewhere. 
For the classification of cf-$m$-manifolds $(m\leq 3)$ 
we refer to \cite{wolf}. 
We shall use the fact that there are $10$ affine diffeomorphism classes of 
closed flat $3$-manifolds $X_1,...,X_{10}$ 
(see e.g. \cite[Theorems 3.5.5 and 3.5.9]{wolf}] 
for the description of them). 
For the classification of closed $4$-manifolds we refer to 
\cite{bullow} or \cite{hillman}. 
\begin{theorem}\label{thm-top-class-dim4}
There are $14$ diffeomorphism classes of orientable, noncompact 
cf-$4$-manifolds. They are represented by: 
$$ \rrr ^4 , \quad 
S^1 \times \rrr ^3 , \quad 
T^2 \times \rrr ^2 , \quad 
TK, \quad 
and  \quad \Lambda ^3 X_j, \; j=1,...,10 .$$ 
\end{theorem} 

Here $K$ is the Klein bottle and $TK$ is the tangent bundle of $K.$ 
To deal with affine classification of cf-$4$-manifolds 
we need some definitions. 
Consider the action of $\sl (2,\zzz )$ on $T^2$, 
induced by the standard action of  $\gl (2,\zzz )$ on $\rrr ^2$, 
and the arising orbit space $T^2/ \gl (2,\zzz )$. 
Let $\rho $ be the equivalence relation in $T^2$ such that 
$$ (e^{i\alpha }, e^{i\beta })\rho 
(e^{i\alpha ^\prime }, e^{i\beta ^\prime}) $$ 
if and only if 
$e^{i\alpha }= e^{i \epsilon _1 \alpha ^\prime }$ and 
$e^{i\beta }=e^{i\epsilon _2( \beta ^\prime -k \alpha ^\prime )} $  
for some $(\epsilon _1 ,\epsilon _2) \in \{ -1,1 \}^2 $ 
and some $k\in \zzz .$

\begin{theorem}\label{thm-aff-class-dim4}
Affine equivalence classes of orientable 
noncompact  cf-$4$-manifolds, 
not diffeomorphic to $T^2\times \rrr ^2$, $TK,$ 
or $S^1\times \rrr ^3,$  
are represented by 
$\rrr ^4$ and $\Lambda ^3 X_j, j=1,...,10.$ 
Affine equivalence classes of cf-manifolds, 
diffeomorphic to $T^2\times \rrr ^2$, $TK,$ 
or $S^1\times \rrr ^3$ correspond to the elements of 
$T^2/\gl (2,\zzz )$, 
$T^2/ \rho $, and 
$S^1/z\sim -z .$ 
\end{theorem} 

%%%%%%%%%%%%%%%%%%%%%%%%%%%%%%%%%%%%
\section[Diffeomorphism classes of some cf-manifolds]
{Diffeomorphism classes of some cf-manifolds}

In this chapter we discuss the problem of the 
topological classification of cf-manifolds in a 
more general context than in Chapter \ref{sec:classif-low-dim}. 
We consider cf-manifolds  homotopy equivalent to 
some low-dimensional Bieberbach manifolds. 
We also deal with stable diffeomorphism classes of some cf-manifolds. 

\begin{defn}\label{def-stble-diffeom} 
Two manifolds $M_1$ and $M_2$ are {\it stably diffeomorphic} 
if there is a positive integer $k$ such that 
$M_1\times \rrr ^k$ is diffeomorphic to 
$M_2\times \rrr ^k.$   
\end{defn}

Given a Bieberbach manifold $X$ and $f_0\in {\rm Aut }\, (\pi _1(X)),$ 
there is a diffeomorphism $f \in {\rm Diff }(X) $ 
 such that $f_\ast =f_0 $ 
(\cite[ch. 2, Theorem 5.3]{charlapb}, \cite[ch. 3, Theorem 3.2.2]{wolf}). 
This induces an action of ${\rm Aut }\, (\pi _1(X))$ on 
$\wt KO(X).$ 
Let $[\xi ]^\ast _\rrr $ be the class of 
a flat bundle $\xi $ over $X$ in 
$\tors \wt {KO}(X)/{\rm Aut}(\pi _1(X)).$ 
The investigation of stable diffeomorphism classes of cf-manifolds  
is based on the following consequence of a result of Mazur 
(cf. \cite[Theorem 2]{mazur}). 

\begin{proposition}\label{prop-KO-descr-of-stable-diff}  
Let $M_1,$ $M_2$ be two cf-$m$-manifolds 
homotopy equivalent to the same Bieberbach manifold $X$ 
and let $\xi _j : M_j \to X, j=1,2,$ 
be the arising flat bundles. 
Assume that $m >2\dim X.$ 
Then the following conditions are equivalent 

\vskip0.1cm 
\noindent 
{\bf a)} $M_1$ and $M_2$ are diffeomorphic, 

\vskip0.1cm 
\noindent 
{\bf b)} 
$[\xi _1 ]^\ast _\rrr = [\xi _2]^\ast _\rrr .$ 
\end{proposition} 

As an immediate consequence of Proposition 
\ref{prop-KO-descr-of-stable-diff} we have 

\begin{corollary}\label{cor-KO-and-stable-diff}
Let $M_1,$ $M_2$ be two cf-$m$-manifolds 
homotopy equivalent to the same Bieberbach manifold $X$ 
and let $\xi _j : M_j \to X, j=1,2,$ 
be the arising flat bundles. 
Then the following conditions are equivalent 

\vskip0.1cm 
\noindent 
{\bf a)} $M_1$ and $M_2$ are stably diffeomorphic, 

\vskip0.1cm 
\noindent 
{\bf b)} 
$[\xi _1 ]^\ast _\rrr =[\xi _2]^\ast _\rrr .$ 
\end{corollary} 

%%%%%%%%
\begin{corollary}\label{cor-SW-classes-and-diffem-of-line-bun}
Diffeomorphism classes of complete flat $(n+1)$-manifolds,  
homotopy equivalent to a fixed Bieberbach $n$-manifold $X$, 
correspond to the elements of 
$H^1(X, \zzz _2)/{\rm Aut}(\pi _1(X)).$ 
\end{corollary} 

Let $g:\rrr ^2 \ni (x,y) \to (-x,y+1)$, let $\mu $ 
be the M\"obius bundle: 
$$ E\mu =\rrr ^2 /\la g \ra \ni [x,y] \to [y] \in  \rrr /\zzz =S^1 ,$$ 
let $P_j : S^1\times S^1 \ni (z_1,z_2) \to z_j, j=1,2,$ 
$\mu _j =P^\ast _j \mu $, and let 
$a,b$ be the generators of the deck group of the Klein bottle 
$K$ defined by the formulas 

\vskip0.2cm 
\centerline{ 
$ a(x,y)=(x +1,y) \quad {\rm and} \quad 
b(x,y)=(-x, y + \frac { 1}{ 2} ) .$ 
} 
\vskip0.2cm 

Consider the generators $\alpha, \beta $ of 
$H^1(K,\zzz _2)$, dual to the images of 
$a,b$ in $H_1(K,\zzz _2 ) $, 
and line bundles $\lambda _1, \lambda _2 $ 
such that $\w _1(\lambda _1)=\alpha $ and $\w _1(\lambda _2) =\beta .$ 
We have 

\begin{proposition}\label{prop-class-cf-over-S1}  
Let $M$ be a cf-$m$-manifold homotopy equivalent to $S^1.$ 
If $m\geq 2 $, then $M$ is diffeomorphic to 
$ S^1\times \rrr ^{m-1}$ or $E\mu \times \rrr ^{m-2} .$ 
\end{proposition} 

%%%%%%%%%%%%%%%%%% 
\begin{theorem}\label{thm-class-cf-over-T2}  
{\bf a)} If $m\geq 5 $, then 
the diffeomorphism classes of 
cf-$m$-manifolds homotopy equivalent to $T^2$ 
are represented by 
$ T^2\times \rrr ^{m-2},$ 
$E(\mu _1 \oplus \Theta _\rrr ^{m-3}), $ 
$E(\mu _1 \oplus \mu _2 \oplus  \Theta _\rrr ^{m-4}), $ 
and  
$E(\mu _1 \oplus \mu _2 \oplus \mu _1\mu _2 \oplus \Theta _\rrr ^{m-5}). 
$ 

\vskip0.1cm 
\noindent 
{\bf b)} 
Diffeomorphism classes of 
cf-$4$-manifolds homotopy equivalent to $T^2$ 
are represented by 
$ T^2\times \rrr ^{2},$ 
$E(\mu _1 \oplus \Theta _\rrr ^{1}), $ and 
$ E(\mu _1 \oplus \mu _2 ). $ 
\end{theorem} 
%%%%%%%
%%% 
\begin{theorem}\label{thm-class-cf-over-K}  
{\bf a)} If $m\geq 5 $, then 
the diffeomorphism classes of 
cf-$m$-manifolds homotopy equivalent to the Klein bottle $K$  
are represented by 
$ K\times \rrr ^{m-2},$ 
$E(\lambda _1 \oplus \Theta _\rrr ^{m-3}), $ 
$E(\lambda _2 \oplus \Theta _\rrr ^{m-3}, $ 
$ E(\lambda _1 \oplus \lambda _2 \oplus  \Theta _\rrr ^{m-4}), $ 
and\\   
$E(\lambda _2 \oplus \lambda _1\lambda _2 \oplus \Theta _\rrr ^{m-5}). 
$

\noindent 
{\bf b)} 
Diffeomorphism classes of 
cf-$4$-manifolds homotopy equivalent to the Klein bottle 
are represented by 
$ K\times \rrr ^{2}, $ 
$ E(\lambda _1 \oplus \Theta _\rrr ^{1}), $ 
$ E(\lambda _2 \oplus \Theta _\rrr ^{1} ) , $ 
$ E(\lambda _1 \oplus \lambda _2 ), $ 
and  
$ E(\lambda _2 \oplus \lambda _1\lambda _2 ). $ 
\end{theorem} 

%%%%%%%%%%%%%%%%%%%
The proof of Proposition \ref{prop-class-cf-over-S1} is easy. 
To describe the idea of the proofs of the other results 
denote  $T^2$ and $K$ by $X.$   
The $3$ and $4$-dimensional case follows from a direct argument. 
By Lemma \ref{prop-class-cf-over-S1}, the isomorphism classes 
of flat vector bundles 
over $X$ whose dimension is greater than $2$ 
correspond to 
their images in ${\rm Tors}\, KO(X).$ 
Using the Atiyah-Hirzebruch $KO$-spectral sequence 
of the fibration $X\to S^1$ one can check 
that the map 
 $$W : {\rm Tors}\, KO(X) \to H^1(X,\zz 2)\oplus H^2(X,\zz 2)  
\cong \zz 2^3, $$
carrying the class of the bundle $\xi $ onto 
$(\w_1(\xi ), \w _2(\xi ) ) ,$  
is a bijection. Now it suffices to find 
the orbit space of the action of $\pi _1(X)$ on 
$ H^1(X,\zz 2)\oplus H^2(X,\zz 2).$ 

%%%%%%%%%%%%%%%%%%%%
\section{Holonomy homomorphisms and geometric invariants} 

In this section we express characteristic classes  
of some flat bundles in terms of their holonomy homomorphisms. 
We also discuss how to calculate cohomology groups 
containing some invariants arising 
in this paper. 
By Corollary \ref{corsplitcommutativeflatbundles}, 
any cf-manifold with abelian vertical holonomy group is 
the total space of the direct sum of 
complex line bundles $L_1,...,L_k $ 
and real line bundles $\lambda _1,...,\lambda _l.$ 
These line bundles are determined by their Chern classes 
and Stiefel-Whitney classes, respectively. 

\begin{lemma}\label{lem-SW-and-holonomy}
Let $ \lambda :E\lambda \to X$ be a real flat line bundle
over a closed flat manifold $X.$
Assume that $\Phi =\Phi _X \times \Phi _U $ is the holonomy
homomorphism of $\lambda $ and
$ \Phi _X \times \Phi _U =(\Psi _X \times \Psi _U )\circ \Pi _{ab} .$
Let $P:H_1(X,\zzz ) \to H_1(X,\zz 2)$ be the projection and let
$\mu :O(1)\to \zz 2$ be the isomorphism.
Then
\vskip0.17cm 
\begin{center} 
$ w_1(\lambda )\circ P = \mu \circ \Psi _U . $
\end{center} 
\end{lemma}

%%%%%%%%%%%%%%%%%%%%%%%%%%
To state an analogous description of the Chern classes 
write the first homology group of a Bieberbach manifold 
$X=\widetilde X/\Gamma $ 
as $\zzz ^{b_1(X)} \oplus S,$ 
where $S={\rm Tors }\, H_1(X,\zzz ).$ 
Let 
 $k_X$ be the order of the holonomy group of $X$ 
and  
let $S_2$ be the torsion subgroup of $H^2(X, \zzz ).$ 
The set $\flbx $,  of isomorphism classes of flat complex
line bundles over $X$, is a commutative group
with tensor product as a group operation.
It is known that $c_1 : \flbx \to H^2(X,\zzz )$ is a monomorphism
(cf. \cite[ch. 16, Theorem 3.4]{husemoller})
and $c_1(\flbx )=S_2 $  (\cite[Theorem 6.1]{kt}).
%%%%%%%%%%%%%%%
For every 
$\Psi \in {\rm Hom \, } (S, \zz {k_X}) $ the formula
$$ \Psi ^H (x) =\begin{cases}
 \Psi (x) & \text{for} \; x\in S \\
0 & \text{for} \; x\in \zzz ^{b_1(X)} 
\end{cases}
 $$
defines
$\Psi  ^H \in {\rm Hom \, } (H_1(X,\zzz ) , \zz {k_X})=
H^1(X, \zz {k_X}).$
%%%
%%%
Consider the  coboundary homomorphism
$\delta :H^1(X, \zz {k_X}) \to H^2(X, \zzz ) $
induced by the short exact sequence
$$ 0 \to \zzz \buildrel {\lambda } \over \rightarrow
\zzz \buildrel {j } \over \rightarrow \zz {k_X} \to 0
,  \leqno {(*)}$$
where $\lambda (x)=k_Xx $ and $j$ is the canonical projection.
%%%%%%%%%%%%
%%%%%%%%%%%%%%
\begin{lemma}\label{lem-first-chern-and-holonomy}
Let $X$ be a Bieberbach manifold and  
let $k_X,$  $\Psi ^H$, and $\delta $ be as above.
Let $L$ be a complex flat line bundle over $X$ with holonomy 
homomorphism $\Phi _L$. 
Take the factorization 
$$ \Phi _L :\Gamma \buildrel {\Pi _{ab}   } \over {\longrightarrow }
H_1(\Gamma ) \buildrel {\Psi _{L} } \over {\longrightarrow }
{\rm U}(1)  $$ 
of $\Phi _L$. 
Then 

\vskip0.1cm 
\noindent 
{\bf a)} $ \delta _S : {\rm Hom }\, (S, \zz {k_X}) \to \tors \, H^2(X,\zzz )$
carrying 
$ \Psi \in {\rm Hom }\, (S, \zz {k_X})$ 
onto $ \delta (\Psi ^H )$ is an isomorphism, 

\vskip0.1cm 
\noindent 
{\bf b)} 
$ c_1(L ) = \delta _S (\Psi _L\vert _S  )  . $
\end{lemma}

The proof of Lemma \ref{lem-SW-and-holonomy} is an easy exercise. 
We do not know a reference to the statement and proof of 
Lemma \ref{lem-first-chern-and-holonomy}.  

%%%%%%%%%%% 
%%%%%%%%%%%%%%%%%%
\begin{corollary}\label{cor-first-chern-and-tor-H1} 
{\bf a)} $c_1(\flbx ) \cong {\rm Hom \, }(S,\zz {k_X}).$ 

\vskip0.1cm \noindent 
{\bf b)} If $\tors H_1(X,\zzz ) =\{ 0\} ,$ 
then all complex flat line bundles over $X$ are trivial. 

\vskip0.1cm \noindent 
{\bf c)} If $\tors H_1(X,\zzz ) =\{ 0\} $ and $M$ is a cf-manifold 
homotopy equivalent to $X$, having abelian vertical holonomy group, 
then 
$M$ is diffeomorphic to the total space 
of the direct sum of some real line bundles over $X.$ 
\end{corollary} 

There are different methods allowing to calculate first and second 
cohomology groups of a Bieberbach $n$-manifold $X$. 
One can use a general approach based on the Smith normal form 
of an integer matrix (cf. \cite{kamroz}, \cite{sims}). 
We discuss another simple approach that can be applied 
if the holonomy group of $X$ is a cyclic group of order 
$k$. In this case $X$ is affinely diffeomorphic to the mapping 
torus $M(g)$ of an affine diffeomorphism 
$g:T^{n-1} \to T^{n_1}$ such that $g^k=id.$ 

%%%%%%%%%%%%%%%
\begin{theorem}\label{thm-H1(G,A)-is-TorsH1(G,Z)}
Let $G$ be a group isomorphic to $\zzz _k $ 
acting on $A\cong \zzz ^{n-1} $ 
  and let $g_0$ be a generator of $G.$
Assume that $g$
is an affine diffeomorphism of $T^{n-1} $ such that $g_\ast =g_0$ 
 and  $M(g)$ is the mapping torus of $g.$
Then
$$ \tors H_1(M(g),\zzz ) \cong H^1(G,A)  \cong \tors  A_G .$$
\end{theorem}
%%%%%%%%%%%

Here $A_G=A/{\rm im}(g_0-id).$ 
The cohomology group $H^1(G,A)$ is of interest in its own right. 
To see this denote by $n_{CC}(G,V)$ the number of connected 
components of the fixed point set of a smooth action of 
$G$ on a manifold $V$ homotopy equivalent to $T^{n-1}.$   
Identifying $\pi _1(T^{n-1})$ with $A$ we can 
treat $A$ as a $\zzz [G]$-module. 
If $G$ is a $p$-group, then 
$$n_{CC}(G,V) =\card H^1(G,A)$$ 
(\cite[Theorem A.10]{cora}). 
Recall that a {\it $G $-lattice} 
is a $G$-module that is also 
a free abelian group of finite rank.

%%%%%%%%%%%%%%%%%%%%%%%%%%%%%%%%%%%%%%%%%%%
\begin{theorem}\label{thm-cal-1st-homology}
Let $G$ be a finite group and let $A$ be a $G$-lattice.
Assume that $m = \vert G \vert $
and $q$ is a positive integer relatively prime to $m.$
Then
$$ \card  H^1(G, A) =
\card  (A\otimes \zz m )^G \, 
\, m^{ - {\rm rank} _{\zz q } (A\otimes \zz q)^G }
. $$
\end{theorem}

\begin{corollary}\label{cor-h1-for-zp}
Let $p$ and $q$ be two different prime numbers and let
$A$ be a $\zz p$-lattice. Then
$$ \card  H^1(G, A) =
p^{\dim _{\zz p} (A\otimes \zz p)^{\zz p}  -
\dim _{\zz q } (A\otimes \zz q)^{\zz p}  }
. $$
\end{corollary}

Corollary \ref{cor-h1-for-zp} is particularly convenient because 
it reduces the calculation of $\card  H^1(G, A) $ 
to the determination of the number of 
solutions of  systems of linear equations 
in finite fields $\zzz _p$ and $\zzz _q .$ 
Applying Theorem \ref{thm-H1(G,A)-is-TorsH1(G,Z)}
and Lemma \ref{lem-first-chern-and-holonomy} we have 

\begin{corollary}\label{cor-relat-homol-coh} 
Let $X$ be  closed flat manifold, 
whose holonomy group is isomorphic to $\zzz _{k} ,$ 
let $g$ be an affine diffeomorphism of $T^{n-1}$ such that 
$X$ is diffeomorphic to $M(g)$, 
and let $A$ denote $\pi _1(T^{n-1})$ 
with the induced $\zzz _k$-action on it.  
Then 
$$ {\rm Tors }\, H^2(X,\zzz ) \cong 
H^1(\zzz _k,A) \quad {\it and } \quad 
 H^1(X,\zzz _2) \cong \zzz _2^{b_1(X)} \oplus 
H^1(\zzz _k,A)\otimes \zzz _2 .$$ 
\end{corollary}

Now we give a convenient criterion of the triviality of
$H^1(G,A).$

\begin{proposition}\label{prop-conditions-H1-trivial}
Let $G$ be a finite group, let $A$ be a
$G$-lattice, and let $q$ be a prime number such that 
$(\vert G \vert , q ) =1 $ .
Assume that for every prime divisor $p$ of $\vert G \vert $ we have  
$$ \dim _{\zz p }(A\otimes \zz p)^G=
\dim _{\zz q }(A\otimes \zz q)^G. $$ 
Then
$H^1(G,A)=\{ 0\} $. 
\end{proposition}

%%%%%%%%%%%%%%%%%%%%%%%%%%%%%%%%%%%%%%%%%%%%%%
%%%%%%%%%%%%%%%%%%%%%%%%%%%%%%%%%%%%%%%%%%%%%%%%%%%%%%%%%%%
%%%%%%%%%%%             BIBLIOGRAPHY         %%%%%%%%%%%%%%
%%%%%%%%%%%%%%%%%%%%%%%%%%%%%%%%%%%%%%%%%%%%%%%%%%%%%%%%%%%

\vskip0.2cm 


\begin{thebibliography}{21}

\bibitem{bullow}
H. Brown, R. B\"ullow, J. Neb\"user, H. Wondratschek, H.Zassenhaus, 
{\it Crystalographic groups of four dimensional space,}
Wiley, New York, 1978.

\bibitem{charlap} L. Charlap,
{\it Compact flat Riemannian manifolds I,}
Ann. Math. 81(1965), 15-30.

\bibitem{charlapb} L. Charlap,
{\it Bieberbach groups and flat  manifolds,}
Springer-Verlag, New York, 1986.

\bibitem{cg}  J. Cheeger, D. Gromoll,
{\it On the structure of complete manifolds of nonnegative curvature,}
Ann. of Math., 96(1972), 413-443.

 
\bibitem{cidschulz} C. Cid, T. Schulz, Computation of five- and six-dimensional
Bieberbach groups, Experiment. Math. {\bf 10 (1)} 2001, 109-115.

\bibitem{cora} P.E.  Conner, F. Raymond, 
{\it  Manifolds with few periodic homeomorphism},
Proc. 2nd Conf. on Compact Transformation Groups, 
Lecture Notes in Math. 299 (Springer, Berlin 1972) 1-75.

\bibitem{ds}   P. Deligne, D. Sullivan,
{\it Fibr\' es vectoriels complexes \`a groupe structural discret,}
C. R. Acad. Sc. Paris, 281(1975), 1081-1083.

\bibitem{dyer} E. Dyer, {\it  Cohomology theories},
Benjamin, New York 1969.

\bibitem{gh}  W.M. Goldman, M.W. Hirsch,
{\it Flat bundles with solvable holonomy,}
Proc. Amer. Math. Soc. 82(1981), 491-494.

\bibitem{gw}   L. Guijarro, G. Walschap,
{\it The metric projection onto the soul,}
Trans. Amer. Math. Soc. 352(2000), 55-69.

 
\bibitem{hillman} J.A. Hillman, Flat $4$-manifold groups,  
New Zealand J. Math., 
{\bf 21} (1995), 29-40.

\bibitem{ht}  M. Hirsch, W. Thurston,
{\it Foliated bundles, invariant measures and flat nanifolds,}
Ann. of Math. \textbf{101} (2) (1975), 369-390.

\bibitem{husemoller}   D. Husemoller,{\it Fibre bundles,}
Mc-Graw-Hill, New York 1966.

\bibitem{kamroz} T. Kaczynski, K. Mischaikow, M. Mrozek, 
{\it Computational homology}, Springer, Berlin, 2004. 


\bibitem{kt}   F. Kamber, Ph. Tondeur,
{\it Flat bundles and characteristic classes of group-representations,}
Amer. J. Math.,  89 (1967), 857-886.

\bibitem{ktflat}   F. Kamber, Ph. Tondeur,
{\it Flat manifolds,} Lecture Notes in
Math. 67, Springer, Berlin 1968,

\bibitem{karoubi}
M. Karoubi,
 {\it K-theory. An Introduction},
Springer, Berlin 1978.

\bibitem{kn} S. Kobayashi, K. Nomizu,
{\it Foundations of differential geometry,}
Interscience Publishers, v. 1; 1963, v. 2; 1969.

\bibitem{mazur}
 B. Mazur,  {\it Stable equivalence of differentiable manifolds},
Bull. Amer. Math. Soc. 67 (1961), 377-384.

\bibitem{milnor}   J. Milnor, {\it On the existence of a connection
with curvature $0$,} Comment. Math. Helv. 32(1958), 215-223.

\bibitem{milnorstasheff} J. Milnor, J. Stasheff,
{\it Characteristic classes,}
Princeton University Press, Princeton 1974.

\bibitem{ow}  M. \"Ozaydin, G. Walschap, {\it Vector bundles with no soul,}
Proc. Amer. Math. Soc. 120(1994), 565-567.

\bibitem{msb} M. Sadowski,
{\it Affinely equivalent complete flat manifolds},
 Cent. Eur. J. Math. {\bf 2}(2) 2004, 332-338. 

\bibitem{sims} Ch.C. Sims, 
{\it Computation with finitely presented groups},
Cambridge University Press, Cambridge 1994.

\bibitem{tapp-thesis} K. Tapp, 
{\it The geometry of open manifolds of nonnegative 
curvature}, Ph.D. Thesis, University of Pennsylvania, 1999. 

\bibitem{tapp} K. Tapp, 
{\it Finiteness theorems for bundles}, preprint, 
University of Pennsylvania, 2004. 

\bibitem{th}   J.A. Thorpe, {\it Parallelizability and flat manifolds,}
Proc. Amer. Math. Soc. 16(1965), 138-142.

 
\bibitem{wilking} B. Wilking, 
{\it On fundamental groups  of manifolds of nonegative curvature,}  
Differential Geom. Appl. 
{\bf 13} (2) (2000), 129-165. 

\bibitem{wolf}   J. Wolf,
{\it Spaces of constant curvature,}
McGraw-Hill, 1967.
\end{thebibliography}
\end{document}